\documentclass[a4paper]{article}
\usepackage[UKenglish]{babel}
\usepackage[T1]{fontenc}
\usepackage{libertine}
\usepackage{libertinust1math}
\usepackage{inconsolata}
\usepackage{AlegreyaSans}
\usepackage{titlesec}
\titleformat{\section}{\normalfont\sffamily\Large\bfseries}{\thesection}{1em}{}
\titleformat{\subsection}{\normalfont\sffamily\large\bfseries}{\thesubsection}{1em}{}

\usepackage[colorlinks,linkcolor=blue,citecolor=black,urlcolor=cyan]{hyperref}
\usepackage{cleveref}
\usepackage{microtype}
\usepackage{fancyvrb}
\usepackage{abstract}

\usepackage{amsmath}
\usepackage{amssymb}
\usepackage{mathtools}
\usepackage{amsthm}
\usepackage{thmtools,thm-restate}
\usepackage{braket}
\usepackage{mathrsfs}
\usepackage{mathbbol}

\usepackage{tikz-cd}

\declaretheorem[numberwithin=section]{theorem}
\declaretheorem[sibling=theorem]{lemma,proposition,corollary}
\declaretheorem[numbered=no,name=\Cref{thm:vbun}]{main-theorem}

\declaretheorem[numbered=no,name=\Cref{prop:flat}]{flatness}
\declaretheorem[numbered=no,name=\Cref{prop:ci}]{completeintersection}
\declaretheorem[numbered=no,name=\Cref{cor:int}]{integral}
\declaretheorem[sibling=theorem,style=definition]{definition}
\declaretheorem[numbered=no,style=remark]{remark}

\newcommand{\A}{\mathbb{A}}
\newcommand{\C}{\mathbb{C}}
\newcommand{\sC}{\mathscr{C}}
\newcommand{\E}{\mathscr{E}}

\newcommand{\h}{\mathrm{H}}
\newcommand{\kk}{\mathbf{k}}
\newcommand{\N}{\mathbb{N}}
\newcommand{\Os}{\mathscr{O}}

\newcommand{\U}{\mathscr{U}}
\newcommand{\Z}{\mathscr{Z}}

\newcommand{\philb}{\hilb^n_{\mathbb{A}^2/\kk}}

\DeclareMathOperator{\hilb}{Hilb}
\DeclareMathOperator{\Hom}{Hom}

\DeclareMathOperator{\rk}{rk}
\DeclareMathOperator{\spec}{Spec}
\DeclareMathOperator{\Spec}{\underline{Spec}}
\DeclareMathOperator{\sym}{Sym}

\title{Hilbert schemes of points on plane curves via standard families}
\author{Miika Rankaviita}
\begin{document}
\maketitle
\begin{abstract}
  We describe an approach to computing
  Hilbert schemes of points on embedded families of plane curves.
  This is achieved by introducing the notion of a standard family,
  and the associated Hilbert schemes are shown to be
  algebraic vector bundles over the Hilbert scheme of points on the plane.
  Any embedded family of plane curves lies in some standard family,
  which allows for the computation of the corresponding Hilbert scheme from that
  of the standard family. Further results on the geometry of these schemes
  are also collected.
\end{abstract}

\section{Introduction}\label{sec:intro}
Throughout, let $\kk$ be a fixed algebraically closed field
of characteristic 0. All schemes are assumed to be defined over $\kk$.
The objects of study in this article are embedded families $\sC\to B$
of affine plane curves and their Hilbert schemes. Here the base $B$
is assumed to be affine and $\sC$ is taken to be a hypersurface in $B\times\A^2$
defined by some regular function
\[f=\sum_{i,j\in\N}f_{ij}\otimes x^iy^j\in\Os_B\otimes_{\kk}\kk[x,y],\quad f_{ij}\in\Os_B.\]
After fixing an $n\in\N$, one can consider the relative Hilbert scheme $\sC^{[n]}\to B$
as defined by Grothendieck \cite{grothendieck}. For a geometric point $z\in B$,
the fibre $\sC^{[n]}_z$ of this family parametrises the length $n$ subschemes of
the corresponding curve $\sC_z$ in the original family.

Hilbert schemes of points on curves are fundamental to the study of curves;
In particular, they are ubiquitous in enumerative geometry
\cite{pt,kst,goettsche-shende}.
Interestingly, their homology is also connected to link invariants
\cite{shende,maulik}.
The easiest case is the Hilbert scheme $C^{[n]}$ of a smooth plane curve $C$.
In this case $C^{[n]}$ is isomorphic to the symmetric power $C^{(n)}$, which
is smooth of dimension $n$.
The theory of Hilbert schemes on reduced curves is developed in \cite{aik,bgs},
where it is shown that if $\sC\to B$ is a family of reduced plane curves,
then $\sC^{[n]}$ is reduced and \textsc{lci} of dimension $n+\dim B$.
If the curves are furthermore integral, then so is the Hilbert scheme.
Moreover, explicit equations for $C^{[n]}$ in $\philb$ are given in \cite{bgs},
when $C$ is a reduced plane curve.

Nodal curves are studied more closely in \cite{ran}, where equations are given for
the Hilbert schemes on the node $\{xy=0\}$. The equations are obtained from those
of the Hilbert schemes on the versal deformation $\{xy=t\}$.
The case of cuspidal curves is done in \cite{lee} following \cite{ran}.
The Hilbert schemes on the non-reduced lines $\{y^\beta=0\}$ are
computed in \cite{luan}, and this is used to describe the irreducible components
of Hilbert schemes on plane curves in general.

The main result of this paper is a construction that allows for computing
relative Hilbert schemes, generalising the results of \cite{ran} to all
embedded families of plane curves as in \Cref{defin:emb}.
The versal deformations are replaced by \emph{standard families} as in \Cref{defin:std}.
Any embedded family $\sC\to B$ of plane curves is the pullback of a suitable
standard family $\sC_I\to B_I$ along some $B\to B_I$, which implies that
the Hilbert scheme $\sC^{[n]}$ is the pullback of $\sC_I^{[n]}$ along
the same $B\to B_I$. Then, \Cref{sec:vbun} describes an explicit procedure
for constructing a locally free sheaf $\E$ on $\philb$ satisfying the following key result.
\begin{main-theorem}
  For a suitable standard family $\sC_I\to B_I$, the relative Hilbert scheme
  $\sC_I^{[n]}$ is isomorphic to the algebraic vector bundle $\Spec_{\philb}(\sym\E)$.
\end{main-theorem}
This provides an explicit description for the family $\sC^{[n]}\to B$
via pulling back along $B\to B_I$.
A similar construction is mentioned in \cite{kst} in the projective setting.
This paper provides more details in the affine case and presents example
computations in \Cref{sec:ex} using \textsc{Singular} \cite{singular}.

It is well-known that the geometry of the Hilbert schemes in higher dimensions
can exhibit arbitrarily bad singularities \cite{murphy} and there are
still open problems regarding the geometry of $\hilb_{\A^3/\kk}$, including
questions about its components and singularities \cite{jelisiejew}.
On the other hand, \Cref{thm:vbun} demonstrates that the geometry of relative
Hilbert schemes on plane curves is quite accessible. \Cref{sec:app} explores
the geometry of the Hilbert schemes in this case, and proves for example
the following results.
\begin{flatness}
  For any embedded family $\sC\to B$ of plane curves, the Hilbert scheme $\sC^{[n]}$
  is flat over $B$.
\end{flatness}
\begin{completeintersection}
  Fix an embedded family $\sC\to B$ of plane curves and $n\in\N$.
  Denote by $\sC_I\to B_I$ the associated standard family.
  If the map $B\to B_I$ is an embedding of a complete intersection,
  then $\sC^{[n]}$ is a complete intersection in $\sC^{[n]}_I$.
\end{completeintersection}
\begin{integral}
  Suppose $\sC\to B$ is an embedded family of plane curves.
  Assume $B$ is integral and denote by $B_{\mathrm{int}}\subseteq B$
  the locus of integral curves. If $B_{\mathrm{int}}$ is open in $B$,
  then $\sC^{[n]}$ is integral.
\end{integral}
This is hardly everything that could be said about the relative
Hilbert schemes on curves, and one can still propose further problems
with attainable solutions.

\subsection*{Acknowledgements}
I thank my supervisor Oscar Kivinen for his support throughout this project
and Ikshu Neithalath for giving thorough and invaluable comments on the final draft.
I would like to also thank the organisers of \emph{Moduli and tensors Warsaw 26-30 October 2026}
for giving the opportunity to present this research at the conference, and also for
the useful discussions. In particular, I thank Joachim Jelisiejew for pointing out that
I should consider the flatness of the Hilbert scheme, which was a crucial point in
\Cref{sec:app}. I also thank Jakub Koncki for sharing his insights about
the Hilbert schemes on curves. The research is supported by a Väisälä project grant of
the Finnish Academy of Science and Letters.

\subsection*{\textsc{Ai} use disclosure}
No \textsc{ai} tools were used at any point of the research process or writing.

\section{Standard families}
The relative Hilbert scheme $\sC^{[n]}$ on an embedded family $\sC\to B$ of plane curves
can be computed from the Hilbert scheme $\sC^{[n]}_I$ on a suitable \emph{standard family}
$\sC_I\to B_I$.
\begin{definition}\label{defin:emb}
  A morphism $\sC\to B$ of affine schemes is said to be an embedded family of plane curves,
  when there is a factorisation
  \[\begin{tikzcd}
      \sC\arrow[hook,r]\dar&B\times\A^2\arrow[dl,"\mathrm{pr}_1"]\\B&
    \end{tikzcd}\]
  such that the top arrow is the embedding of a hypersurface cut out by
  a single regular function
  \[f=\sum_{i,j\in\N}f_{ij}\otimes x^iy^j\in\Os_B\otimes_\kk\kk[x,y],\quad f_{ij}\in\Os_B,\]
  and furthermore, the fibres of $\sC\to B$ are 1-dimensional.
\end{definition}
\begin{definition}\label{defin:std}
  The standard family of a finite index set $I\subset\N^2$ is the family
  with base
  \[B_I=\spec\kk[b_{ij}\mid (i,j)\in I]\]
  and total space
  \[\sC_I=\spec\frac{\kk[b_{ij}\mid (i,j)\in I][x,y]}{\left(\sum_{(i,j)\in I}b_{ij}x^iy^j\right)}.\]
  The structure map $\sC_I\to B_I$ is defined by the obvious inclusion
  of structure rings.
\end{definition}
\begin{remark}
  A standard family defines an embedded family of plane curves
  away from the line in $B_I$ parametrising the constant polynomials.
  The fibre over $0$ is the entire plane, while the fibres over
  non-zero constants are empty. However, when these base points
  are preserved, the projection $\sC_I\to\A^2$ is a vector bundle.
  Thus, the base points are kept for convenience, and this does
  not cause any problems if enough care is taken.
\end{remark}
\begin{definition}
  Given an embedded family $\sC\to B$ as in \Cref{defin:emb}
  and a fixed $n\in\N$, the standard family associated to
  the data of $\sC\to B$ and $n$ is the standard family $\sC_I\to B_I$
  of the index set
  \[I=\Set{(i,j)\in\N^2\mid \text{$i+j<n$ or $f_{ij}\neq 0$}},\]
  where the $f_{ij}$ are the coefficients of the regular function
  defining $\sC\hookrightarrow B\times\A^2$.
\end{definition}
\begin{remark}
  The reason why one would want to also include all those tuples $(i,j)$ with
  $i+j< n$ will become apparent in \Cref{sec:vbun}.
\end{remark}

Next, fix an embedded family $\sC\to B$ of plane curves along
with $n\in\N$, and let $\sC_I\to B_I$ be the standard family
associated to this data. One can then define a map
\[\Os_{B_I}\to\Os_B:b_{ij}\mapsto f_{ij}.\]
This map gives rise to a pullback square
\[\begin{tikzcd}
    \sC\rar\dar\arrow[dr,phantom,very near start,"\lrcorner"] & \sC_I\dar \\
    B\rar & B_I
  \end{tikzcd}\]
Furthermore, this induces a pullback square
\[\begin{tikzcd}
    \sC^{[n]}\rar\dar\arrow[dr,phantom,very near start,"\lrcorner"] & \sC_I^{[n]}\dar \\
    B\rar & B_I
  \end{tikzcd}\]
expressing the Hilbert scheme $\sC^{[n]}$ in terms of that of $\sC_I$.
The next section is devoted to proving that the Hilbert scheme $\sC^{[n]}_I$
is smooth and can be constructed explicitly as an algebraic vector bundle
over $\philb$. This makes it possible to compute $\sC^{[n]}$
by first computing $\sC^{[n]}_I$.

\section{Hilbert schemes of standard families}\label{sec:vbun}
The purpose of this section is to describe the Hilbert scheme $\sC_I^{[n]}$
on a standard family as an algebraic vector bundle over $\philb$.
Denote in this section
\[\h\vcentcolon=\philb\]
for short. The proof rests on the fact that $\sC_I^{[n]}$ is smooth,
which follows easily from the description of the tangent spaces
in \cite[Section~4]{shende}. The argument is reproduced below for
the convenience of the reader.
\begin{lemma}\label{lem:smth}
  Fix a geometric point $\kappa\in\sC^{[n]}_I$. If for some subset $J\subseteq I$,
  the set $\set{x^iy^j\mid(i,j)\in J}$ defines a basis for $\Os_\kappa$,
  then the Hilbert scheme $\sC^{[n]}_I$ is smooth at $\kappa$.
\end{lemma}
\begin{proof}
  As the non-empty fibres of $\sC_I^{[n]}$ have dimension at least $n$,
  one obtains the bound
  \[\dim\sC_I^{[n]}\geq\dim B_I+n.\]
  To conclude that $\sC_I^{[n]}$ is smooth, it suffices to show that
  the tangent spaces have dimension $\dim B_I +n$.

  Fix a geometric point $\kappa\in\sC_I^{[n]}$ such that $\kappa$ defines
  a subscheme of the fibre $\sC_{I,z}$ above the point $z\in B_I$.
  Consider the exact sequence
  \[\begin{tikzcd}
      0\rar&T_\kappa\sC_I^{[n]}\rar&T_zB_I\times T_\kappa\philb\rar
      &\Os_\kappa
    \end{tikzcd}\]
  Here $T_\kappa\philb$ is identified with
  $\Hom_{\kk[x,y]}(I_\kappa,\Os_\kappa)$. The last arrow maps
  a pair $(v,\eta)$ to the reduction of the polynomial
  \[\eta\left(\sum_{(i,j)\in I}z_{ij}x^iy^j\right)-\sum_{(i,j)\in I}v_{i,j}x^iy^j\]
  in $\Os_\kappa$. Since $I$ contains a subset $J$ defining a basis for $\Os_\kappa$,
  the last map is surjective. Thus, counting dimensions gives
  \[
    0=\dim T_\kappa\sC_I^{[n]}-(\dim B_I+2n)+n
    \implies\dim T_\kappa\sC_I^{[n]}=\dim B_I+n.
  \]
\end{proof}

Now, fix an $n\in\N$ and a finite index set $I\subset\N^2$ containing all
the tuples $(i,j)$ such that $i+j< n$. The Hilbert scheme $\sC_I^{[n]}$ on
the standard family can be constructed as an algebraic vector bundle over
$\h=\philb$. The main observation is that for any geometric point $\kappa\in\h$,
it is possible to solve for all the $z\in B_I$ such that $\kappa$ defines
a subscheme of the fibre $\sC_{I,z}$ above $z\in B_I$.
Indeed, recall that $\sC_{I,z}$ is the vanishing locus of the polynomial
$\sum_{(i,j)\in I} z_{ij}x^iy^j$ in $\A^2$, so $\kappa\subset\sC_{I,z}$ if and only if
the polynomial $\sum z_{ij}x^iy^j$ reduces to $0$ in $\Os_\kappa$. Thus,
the kernel of the composition
\[\bigoplus_{(i,j)\in I}\kk\cdot x^iy^j\hookrightarrow\kk[x,y]\to\Os_\kappa.\]
of $\kk$-linear maps parametrises those $z\in B_I$ for which $\kappa\subset\sC_{I,z}$.

This observation is globalised by considering the composition
\begin{equation}\label{eq:rho}\tag{\dag}
\rho\colon\bigoplus_{(i,j)\in I}\Os_{\h}x^iy^j\hookrightarrow\pi_\ast\Os_{\h\times\A^2}
\to\pi_\ast\Os_\U
\end{equation}
of $\Os_\h$-module morphisms, where $\pi$ is the projection $\h\times\A^2\to\h$
and $\U\subset\h\times\A^2$ is the universal family over $\h$. Denote
$\E\vcentcolon=(\ker\rho)^\vee$ and define the algebraic vector bundle
\begin{equation}\label{eq:vbun}\tag{\ddag}
  X\vcentcolon=\Spec_\h(\sym\E)
\end{equation}
over $\h$.

This vector bundle embeds in $\h\times B_I$. Indeed, consider the sequence
\[\begin{tikzcd}
    0\rar&\ker\rho\rar&\displaystyle\bigoplus_{(i,j)\in I}\Os_\h x^iy^j\rar{\rho}
    &\pi_\ast\Os_\U\rar&0
\end{tikzcd}\]
The fibres of $\pi_\ast\Os_\U$ above geometric points of $\h$ are
length $n$ quotients of $\kk[x,y]$, so they are spanned by
the set of monomials $x^iy^j$ with $i+j< n$. Since $I$ contains
these exponents, one can conclude $\rho$ is surjective (by Nakayama's lemma)
and the above sequence is exact. Next, taking duals gives a \textsc{ses}
\[\begin{tikzcd}
    0\rar&\left(\pi_\ast\Os_\U\right)^\vee\rar
    &\left(\bigoplus\Os_\h x^iy^j\right)^\vee\rar
    &\left(\ker\rho\right)^\vee\rar&0
\end{tikzcd}\]
as $\pi_\ast\Os_\U$ is flat. Then, taking $\sym$ and relative $\Spec$
of the second arrow gives a closed embedding
\[\Spec_\h(\sym\E)\hookrightarrow
\Spec_\h\left(\sym\left(\bigoplus\Os_\h x^iy^j\right)^\vee\right)\]
The target is isomorphic to $\h\times B_I$, so this is a closed embedding
$X\hookrightarrow\h\times B_I$. The set of $\kk$-points of $X$ can
be identified as a subset of the set of $\kk$-points of $\h\times B_I$.
\begin{lemma}\label{lem:pts}
  The set of $\kk$-points of $X$ is
  \[X(\kk)=\Set{(\kappa,z)\in\h(\kk)\times B_I(\kk)\mid \kappa\subset\sC_{I,z}}.\]
\end{lemma}
\begin{proof}
  Fix a $\kk$-point $x\colon\spec\kk\to X$, and denote the projections to $\h$
  and $B_I$ by $\kappa\colon\spec\kk\to\h$ and $z\colon\spec\kk\to B_I$ respectively.
  In order to show that $\kappa\subset\sC_{I,z}$, begin by restricting to fibres at $\kappa$.
  In other words, take the following pullbacks.
  \[\begin{tikzcd}[row sep=2.5em]
      X\arrow[r,hookrightarrow]&\h\times B_I\rar&\h\\
      X\vert_\kappa\arrow[u,hookrightarrow]\arrow[r,hookrightarrow]
      \arrow[ur,phantom,very near start,"\urcorner"]&B_I\rar
      \arrow[u,hookrightarrow]\arrow[ur,phantom,very near start,"\urcorner"]
      &\spec\kk\arrow[u,hookrightarrow,"\kappa"']
    \end{tikzcd}\]
  The embedding $X\vert_\kappa\hookrightarrow B_I$ can be understood
  via the following construction. First, restrict the \textsc{ses}
  \[\begin{tikzcd}
      0\rar&\ker\rho\rar&\displaystyle\bigoplus_{(i,j)\in I}\Os_\h x^iy^j\rar{\rho}
      &\pi_\ast\Os_\U\rar&0
    \end{tikzcd}\]
  to fibres at $\kappa$ and take duals, obtaining the \textsc{ses}
  \[\begin{tikzcd}
      0\rar&(\Os_\kappa)^\vee\rar{\rho\vert_\kappa^\vee}&(\bigoplus\kk\cdot x^iy^j)^\vee\rar
      &\E\vert_\kappa\rar&0
    \end{tikzcd}\]
  Exactness follows from the fact that $\pi_\ast\Os_\U$ is flat. The embedding
  $X\vert_\kappa\hookrightarrow B_I$ is the result of taking $\sym$ and $\spec$ of
  the latter arrow.

  Next, observe that $x\colon\spec\kk\to X$ induces a lift
  \[\begin{tikzcd}
      X\vert_\kappa\arrow[r,hookrightarrow]&B_I\\
      \spec\kk\arrow[u,dashrightarrow]\arrow[ur,"z"']
    \end{tikzcd}\]
  This gives rise to a triangle
  \[\begin{tikzcd}
      \E\vert_\kappa\dar&(\bigoplus\kk\cdot x^iy^j)^\vee\lar\arrow[dl,"z^\#"]\\
      \kk&
    \end{tikzcd}\]
  of $\kk$-vector spaces. Now, take an arbitrary $L\in (\bigoplus\kk\cdot x^iy^j)^\vee$
  that maps to $0$ in $\E\vert_\kappa$. Observe that
  \[0=z^\#(L)=L\left(\sum_{(i,j)\in I}z_{ij}x^iy^j\right).\]
  Now, since $L$ was arbitrary and the kernel of
  $(\bigoplus\kk\cdot x^iy^j)^\vee\to \E\vert_\kappa$ can be identified with $(\Os_\kappa)^\vee$,
  this is implies
  \[\rho\vert_\kappa\left(\sum_{(i,j)\in I}z_{ij}x^iy^j\right)=0\]
  in $\Os_\kappa$. Therefore, $\kappa\subset\sC_{I,z}$.

  Conversely, suppose $(\kappa,z)\in\h(\kk)\times B_I(\kk)$ satisfies $\kappa\subset\sC_{I,z}$.
  Then, retracing the steps above backwards shows $z\colon\spec\kk\to B_I$ lifts to $X\vert_\kappa$.
\end{proof}
\begin{theorem}\label{thm:vbun}
  If $I\subset\N^2$ is a finite index set containing all
  the tuples $(i,j)$ such that $i+j< n$, then
  the relative Hilbert scheme $\sC^{[n]}_I$ is isomorphic to
  the algebraic vector bundle $X=\Spec_{\hilb^n_{\A^2/\kk}}(\sym\E)$ as defined at \eqref{eq:vbun} on
  page~\pageref{eq:vbun}.
\end{theorem}
\begin{proof}
  As above, denote $\h\vcentcolon=\philb$. Take the family
  \[\Z\vcentcolon=(X\times_{B_I}\sC_I)\times_{\h\times\A^2}\U\]
  as a candidate $X$-point of $\sC_I^{[n]}$ defining an isomorphism
  $X\to\sC_I^{[n]}$. The following commutative diagram will be useful.
  \[\begin{tikzcd}[row sep=3em]
      \Z\arrow[r,hookrightarrow]\dar\arrow[dr,phantom,very near start,"\lrcorner"]
      &X\times_\h\U\rar\dar\arrow[dr,phantom,very near start,"\lrcorner"]
      &\U\arrow[d,hookrightarrow]\\
      X\times_{B_I}\sC_I\arrow[r,hookrightarrow]\drar
      &X\times\A^2\rar\dar\arrow[dr,phantom,very near start,"\lrcorner"]
      &\h\times\A^2\dar\\&X\rar&\h
    \end{tikzcd}\]
  In order to verify that $\Z$ indeed defines an $X$-point of $\sC_I^{[n]}$,
  one must check that $\Z\to X\times_{B_I}\sC_I$ is a closed embedding,
  the composition $\Z\to X$ is flat and the Hilbert polynomial of
  a geometric fibre is the constant $n$. Since $\Z\to X\times_{B_I}\sC_I$ is
  the pullback of a closed embedding, the first property follows immediately.
  Next, inspect the geometric fibres of the composition
  \[\Z\hookrightarrow X\times_{B_I}\sC_I\to X.\]
  Fix a geometric point $x\colon\spec\kk\to X$, and denote the projections to $\h$
  and $B_I$ by $\kappa\colon\spec\kk\to\h$ and $z\colon\spec\kk\to B_I$ respectively.
  According to the above diagram, the fibre $\Z_x$ is the pullback
  \[\begin{tikzcd}
      \Z_x\rar\dar\arrow[dr,phantom,very near start,"\lrcorner"]
      & \kappa\arrow[d,hookrightarrow] \\ \sC_{I,z}\arrow[r,hookrightarrow] & \A^2
    \end{tikzcd}\]
  By \Cref{lem:pts}, there is a lift
  \[\begin{tikzcd}
      &\kappa\arrow[dl,dashrightarrow]\arrow[d,hookrightarrow]\\
      \sC_{I,z}\arrow[r,hookrightarrow]&\A^2
    \end{tikzcd}\]
  which implies the top arrow $\Z_x\to\kappa$ in the pullback square is an isomorphism.
  Therefore, the fibres of the map $\Z\to X$ have the same Hilbert polynomial as $\kappa$.
  It remains to show that the map is flat. This can be seen from the first commutative diagram.
  The map $\Z\hookrightarrow X\times_\h\U$ is finite as it is a closed embedding.
  The map $X\times_\h\U\to X$ is also finite, since it is the pullback of the finite map $\U\to\h$.
  Thus, the composition $\Z\to X$ is finite. As $X$ is reduced and the Hilbert polynomial
  is constant in the fibres, the map $\Z\to X$ is flat.

  Therefore, $\Z$ defines a map $\Phi\colon X\to\sC_I^{[n]}$.
  Now, $\h$ is smooth and connected \cite{fogarty} and so is $X$. The Hilbert scheme
  $\sC_I^{[n]}$ is also smooth by \Cref{lem:smth}, so by Zariski's Main Theorem
  it remains to note that $\Phi$ is bijective on geometric points.
  Indeed, if one uses \Cref{lem:pts} to identify the $\kk$-points of $X$ with pairs
  $(\kappa,z)\in\h(\kk)\times B_I(\kk)$ such that $\kappa\subset\sC_{I,z}$, then $\Phi$ sends such
  a pair to the $\kk$-point $\kappa\subset\sC_{I,z}$ of $\sC_I^{[n]}$.
\end{proof}

\section{Applications: The geometry of Hilbert schemes}\label{sec:app}
This section collects results about the relative Hilbert schemes on curves
that follow from \Cref{thm:vbun} in combination with
earlier results from the literature. An immediate result is that $\sC^{[n]}$
is flat over the base, which in turn yields further results on
the components of $\sC^{[n]}$.
\begin{proposition}\label{prop:flat}
  For any embedded family $\sC\to B$ of plane curves, the Hilbert scheme $\sC^{[n]}$
  is flat over $B$.
\end{proposition}
\begin{proof}
  Fixing $n\in \N$, any embedded family $\sC\to B$ is the pullback of
  the associated standard family $\sC_I\to B_I$.
  Now, let $B^1_I\subseteq B_I$ be the open locus of parameters corresponding
  to curves in $\sC_I$. Then, the restriction of $\sC_I^{[n]}\to B_I$ to $B_I^1$ is
  a morphism of smooth schemes. Furthermore, the fibre dimension
  is constant above $B_I^1$ \cite{luan}, so the restriction is flat by Miracle flatness.
  The family $\sC^{[n]}\to B$ is the pullback of this restriction.
\end{proof}
One can also arrange the Hilbert scheme to be a complete intersection.
\begin{proposition}\label{prop:ci}
  Fix an embedded family $\sC\to B$ of plane curves and $n\in\N$.
  Denote by $\sC_I\to B_I$ the associated standard family.
  If the map $B\to B_I$ is an embedding of a complete intersection,
  then $\sC^{[n]}$ is a complete intersection in $\sC^{[n]}_I$.
\end{proposition}
\begin{proof}
  Denote $c\vcentcolon=\dim B_I - \dim B$. By assumption, $B$ is cut out
  by regular functions $f_1,\ldots, f_c\in\Os_{B_I}$.
  Then, $\sC^{[n]}$ is cut out by
  \[f_1\otimes1,\ldots,f_c\otimes1\in\Os_{B_I}\otimes\Os_{\philb}.\]
  Now, the fibres of the pullback have dimension $n$ \cite{luan}
  and the base has dimension $\lvert I\rvert - c$. Hence, the total space
  has dimension $\lvert I\rvert +n - c$ so its codimension is $c$.
\end{proof}
The following statement is known, but it is proved for the convenience
of the reader.
\begin{proposition}\label{prop:smth}
  If $\sC\to B$ is an embedded family of plane curves
  over a smooth base, where each fibre is a smooth curve,
  then the Hilbert scheme $\sC^{[n]}$ is smooth.
\end{proposition}
\begin{proof}
  To verify that $\sC^{[n]}$ is smooth at a geometric
  point $\kappa\in\sC^{[n]}$, compute the dimension of
  the tangent space at $\kappa$. Let $z\in B$ such that
  $\kappa\subset\sC_z$, and consider the completion of
  the local ring $\Os_{B,z}$ with maximal ideal $\mathfrak{m}$.
  Denote
  \[
    \sC_D\vcentcolon=
    \spec\left(\frac{\hat{\Os}_{B,z}}{\mathfrak{m}^2}\right)
    \times_{B}\sC.
  \]
  Note that $\sC$ is flat over $B$ by \Cref{prop:flat}
  (alternatively, use the slicing criterion \cite[Theorem~24.6.9]{vakil}).
  Therefore, the morphism
  \[\sC_D\to \spec\left(\hat{\Os}_{B,z}/\mathfrak{m}^2\right)\]
  is also flat, so it can be seen as an infinitesimal deformation of
  the curve $\sC_z$. Since the curve is smooth, the deformation is trivial, implying
  \[
    \sC_D^{[n]}\cong
    \spec\left(\frac{\hat{\Os}_{B,z}}{\mathfrak{m}^2}\right)
    \times\sC_z^{[n]}.
  \]
  Finally, since $\sC_D^{[n]}$ and $\sC^{[n]}$
  have isomorphic tangent spaces at $\kappa$,
  \[
    \dim T_\kappa\sC^{[n]}
    =\dim\left(\frac{\hat{\Os}_{B,z}}{\mathfrak{m}^2}\right)
    +\dim T_\kappa\sC_z^{[n]}
    =\dim B+n.
  \]
\end{proof}

The following results on the components of the Hilbert schemes
can be derived from flatness.
\begin{corollary}\label{cor:dense}
  Suppose $\sC\to B$ is an embedded family of plane curves
  and $B_+\subseteq B$ is a Zariski dense open subscheme. Then,
  the restriction $\sC^{[n]}_+$ of $\sC^{[n]}$ to $B_+$
  is Zariski dense.
\end{corollary}
\begin{proof}
  The question is topological, so one may assume $B$ is reduced.
  Then, as $B_+$ is a dense open, it contains all associated points
  of $B$. Since $\sC^{[n]}\to B$ is flat by \Cref{prop:flat},
  the associated points of $\sC^{[n]}$ map to associated points of $B$.
  Thus, $\sC^{[n]}_+$ contains all associated points of $\sC^{[n]}$.
\end{proof}
\begin{corollary}\label{cor:int}
  Suppose $\sC\to B$ is an embedded family of plane curves.
  Assume $B$ is integral and denote by $B_{\mathrm{int}}\subseteq B$
  the locus of integral curves. If $B_{\mathrm{int}}$ is open in $B$,
  then $\sC^{[n]}$ is integral.
\end{corollary}
\begin{proof}
  Firstly, as $\sC^{[n]}\to B$ is flat by \Cref{prop:flat}
  and $B$ is integral, associated points of $\sC^{[n]}$ map to
  the generic point of $B$. Thus, the irreducible components of $\sC^{[n]}$
  are in bijection with the irreducible components of the generic fibre.
  If the generic fibre admitted more than one irreducible component,
  then there would be a non-empty open subset of $B$, where the fibres are reducible
  \cite[\href{https://stacks.math.columbia.edu/tag/0554}{Lemma 0554}]{stacks-project}.
  Since $B_{\mathrm{int}}$ is a dense open of $B$ and the fibres of $\sC^{[n]}$
  above $B_{\mathrm{int}}$ are integral \cite{aik,bgs}, this would
  lead to a contradiction. Therefore, $\sC^{[n]}$ has a unique irreducible component

  To prove $\sC^{[n]}$ is reduced, it suffices to show it is $R_0+S_1$.
  Begin by noting that a general $\kappa\in\sC^{[n]}$ lies on an integral
  curve and is supported on the smooth locus. Thus, the formal neighbourhood
  of $\kappa$ in $\sC^{[n]}$ is equivalent to that of the Hilbert scheme of points on
  some family of smooth curves. Therefore, $\kappa$ is a smooth point by
  \Cref{prop:smth}, and one can conclude $\sC^{[n]}$ is $R_0$.
  Finally, since $\sC^{[n]}\to B$ is flat and $B$ is reduced,
  $\sC^{[n]}$ cannot contain any embedded points. This implies it is $S_1$
  \cite[\href{https://stacks.math.columbia.edu/tag/031Q}{Lemma 031Q}]{stacks-project}.
\end{proof}

\section{Explicit constructions}\label{sec:constr}
The description of $\sC^{[n]}_I$ as an algebraic vector bundle
over $\philb$ is explicit enough to allow for finding equations describing it as
a subscheme of $\philb\times B$, and also finding a parametrisation
$\philb\times\A^m\xlongrightarrow{\sim}\sC^{[n]}_I$. To obtain explicit equations,
one must restrict to an open cover of $\philb$ that trivialises the universal
family $\U$. Such a cover parametrised by partitions of $n$ is defined in \cite{haiman}.
Fixing a trivialising open $U\subset\philb$, the first proposition below describes the equations of
\[\sC^{[n]}_I\vert_U\vcentcolon=U\times_{\hilb^n_{\A^2/\kk}}\sC^{[n]}_I\]
as a subscheme of $U\times B$.

\begin{proposition}\label{prop:eqns}
  For $n$ and $I$ as in \Cref{thm:vbun} and
  a trivialising open $U\subset\philb$ for $\U$,
  the subscheme $\sC^{[n]}_I\vert_U$ of $U\times B$ is cut out by the equations
  \[
    F_k\vcentcolon=\sum_{(i,j)\in I}\left(\rho\vert_U(x^iy^j)\right)_k\otimes b_{ij}
    \in\Os_U\otimes_\kk\Os_B,\quad 1\leq k\leq n.
  \]
  Here $\rho$ is the morphism defined at \eqref{eq:rho} on page~\pageref{eq:rho},
  and the subscript $(-)_k$ refers to taking the $k$-th component of the section of
  $\pi_\ast\Os_\U\vert_U$ when trivialised by some fixed isomorphism
  $\pi_\ast\Os_\U\vert_U\cong\Os_U^{\oplus n}$.
\end{proposition}
\begin{proof}
  First, consider the \textsc{ses}
  \[\begin{tikzcd}
      0\rar&(\ker\rho)\vert_U\rar
      &\displaystyle\bigoplus_{(i,j)\in I}\Os_Ux^iy^j\rar{\rho}
      &(\pi_\ast\Os_\U)\vert_U\rar&0
    \end{tikzcd}\]
  \Cref{thm:vbun} implies that the map
  \[\sym(\bigoplus\Os_Ux^iy^j)^\vee\twoheadrightarrow\sym(\E\vert_U)\]
  obtained from the above \textsc{ses} by taking duals and $\sym$
  defines the embedding $\sC^{[n]}_I\hookrightarrow U\times B$.
  Now, \cite[\href{https://stacks.math.columbia.edu/tag/00DO}{Lemma 00DO}]{stacks-project}
  implies the ideal of this embedding is generated by the image of
  the map
  \[
    (\pi_\ast\Os_\U)\vert_U^\vee\hookrightarrow\sym(\bigoplus\Os_Ux^iy^j)^\vee
    \xlongrightarrow{\sim}\Os_U\otimes_\kk\Os_B
  \]
  Fix an isomorphism
  \[(\pi_\ast\Os_\U)\vert_U\xlongrightarrow{\sim}\Os_U^{\oplus n}\]
  trivialising $\U$. Denote by $e_k^\vee\in(\Os_U^{\oplus n})^\vee$ the duals of
  the canonical generators $e_1,\ldots,e_n$ of $\Os_U^{\oplus n}$. Then, the ideal of
  the embedding $\sC^{[n]}_I\hookrightarrow U\times B$ is generated by the images of
  the $e_k^\vee$ in $\Os_U\otimes\Os_B$. One can check that
  the image of $e_k^\vee$ in $\Os_U\otimes\Os_B$ is the section
  \[\sum_{(i,j)\in I}\left(\rho\vert_U(x^iy^j)\right)_k\otimes b_{ij}\in\Os_U\otimes\Os_B\]
  as stated.
\end{proof}
For computing a parametrisation $U\times\A^m\to \sC_I^{[n]}\vert_U$, the open set $U$
must trivialise the module $\ker\rho$. This could be obtained by representing
$\rho$ over an open set trivialising $\U$, and finding a trivialising cover for $\ker\rho$
there.
\begin{proposition}\label{prop:param}
  Suppose $n$ and $I$ are as in \Cref{thm:vbun} and set $m\vcentcolon=\lvert I\rvert-n$.
  Let $\rho$ be the morphism defined at \eqref{eq:rho} on page~\pageref{eq:rho}.
  Fix a trivialising open $U\subset\philb$ for $\ker\rho$.
  Take any set of global generators $\sigma^1,\ldots,\sigma^m$ for $(\ker\rho)\vert_U$, where
  \[\sigma^k=\sum_{(i,j)\in I}\sigma^k_{ij}x^iy^j,\quad 1\leq k\leq m\]
  and $\sigma^k_{ij}\in\Os_U$. Then, the homomorphism
  \[\Os_U\otimes_\kk\Os_B\to \Os_U[t_1,\ldots,t_m] : g\otimes b_{ij}\mapsto\sum_{k=1}^mg\sigma^k_{ij}t_k\]
  defines a map $U\times\A^m\to U\times B$ that restricts to an isomorphism
  \[U\times\A^m\xlongrightarrow{\sim}\sC^{[n]}_I\vert_U.\]
\end{proposition}
\begin{proof}
  As above, consider first the \textsc{ses}
  \[\begin{tikzcd}
      0\rar&(\ker\rho)\vert_U\rar
      &\displaystyle\bigoplus_{(i,j)\in I}\Os_Ux^iy^j\rar{\rho}
      &(\pi_\ast\Os_\U)\vert_U\rar&0
    \end{tikzcd}\]
  Counting ranks, gives $\rk(\ker\rho\vert_U)=\lvert I\rvert-n=m$.
  Now, fix an isomorphism
  \[\Os_U^{\oplus m}\xlongrightarrow{\sim}(\ker\rho)\vert_U\]
  sending the $k$-th generator of $\Os_U^{\oplus m}$ to a section $\sigma_k$.
  Since $(\ker\rho)\vert_U$ is a submodule of $\bigoplus\Os_Ux^iy^j$,
  these sections can be written as
  \[\sigma_k=\sum_{(i,j)\in I}\sigma^k_{ij}x^iy^j\]
  with $\sigma^k_{ij}\in\Os_U$.

  Next, taking duals and $\sym$ gives the morphism
  \[\begin{tikzcd}
      \sym(\bigoplus\Os_Ux^iy^j)^\vee\arrow[r,twoheadrightarrow]
      &\sym(\E\vert_U)\rar{\sim}&\sym(\Os_U^{\oplus m})^\vee
    \end{tikzcd}\]
  Using the isomorphism
  $\spec(\sym(\E\vert_U))\cong\sC^{[n]}_I\vert_U$ coming from \Cref{thm:vbun},
  these morphisms define a commutative triangle
  \[\begin{tikzcd}[row sep=3em]
      U\times\A^m\rar\arrow[d,"\sim"']&U\times B\\
      \sC^{[n]}_I\vert_U\arrow[ur,hookrightarrow]&
    \end{tikzcd}\]
  Finally, one can check that the map $U\times\A^m\to U\times B$
  is defined via
  \[\Os_U\otimes_\kk\Os_B\to\Os_U[t_1,\ldots,t_m]:g\otimes b_{ij}\mapsto \sum_{k=1}^mg\sigma^k_{ij}t_k\]
  as stated.
\end{proof}

\section{Examples}\label{sec:ex}
With the results of \Cref{sec:constr}, it is possible to compute
explicit examples.

\subsection{Family of conics}
As a first example, consider the family with base
$B=\spec\C[\alpha,\beta,\gamma]$ and
\[\sC=\spec\frac{\C[\alpha,\beta,\gamma,x,y]}{\alpha+\beta x^2+\gamma y^2}.\]
This is an embedded family of curves away from the locus $\set{\beta=\gamma=0}$,
and a general fibre is a non-degenerate conic. On the locus
$\set{\alpha=0}$ the fibres consist of a pair of lines meeting at the origin.
On the loci $\set{\beta=0}$ and $\set{\gamma=0}$ the fibres are unions
of pairs of lines. All of these loci contain one or both of the non-reduced
axes: $\spec(\C[x,y]/(x^2))$ or $\spec(\C[x,y]/(y^2))$.

Consider the index set $I=\set{(0,0),(1,0),(0,1),(2,0),(0,2)}$.
Then, the Hilbert scheme $\sC^{[2]}$ is the pullback of $\sC_I^{[2]}$
along the obvious inclusion $B\hookrightarrow B_I$. Now, the Hilbert scheme
$\hilb^2_{\A^2/\C}$ is covered by two affine opens corresponding to the partitions of 2.
Denote the partitions by $\lambda=1^2$ and $\mu=2$. Concentrate on the chart $U_\lambda$.
The computations are similar for $U_\mu$. It is easy to check that
\[
  U_\lambda\cong\spec\C[\sigma,\tau,z,c]
\]
and the restriction of the universal family to $U_\lambda$ is defined by the ideal
\[
  (x^2-\sigma x-\tau,y-z x-c)
  \triangleleft\Os_{U_\lambda}\otimes_\C\C[x,y]
\]
Therefore,
\[
  \pi_\ast\Os_\U\vert_{U_\lambda}
  =\frac{\Os_{U_\lambda}\otimes_\C\C[x,y]}{(x^2-\sigma x-\tau, y-z x-c)}
  \cong \Os_{U_\lambda}\oplus\Os_{U_\lambda}x.
\]
It is also straight-forward to verify that the matrix for the morphism
\[
  \rho\colon
  \Os_{U_\lambda}
  \oplus\Os_{U_\lambda}x
  \oplus\Os_{U_\lambda}y
  \oplus\Os_{U_\lambda}x^2
  \oplus\Os_{U_\lambda}y^2
  \to\pi_\ast\Os_\U\vert_{U_\lambda}
\]
is given by
\[\begin{pmatrix}
  1&0&c&\tau&\tau z^2+c^2\\
  0&1&z&\sigma&\sigma z^2+2z c
\end{pmatrix}\]
By \Cref{prop:eqns}, the intersection $\sC^{[2]}_I\cap(U_\lambda\times B_I)$
is cut out by the ideal
\[
  (b_{00}+c b_{01}+\tau b_{20}+(\tau z^2+c^2)b_{02},
  b_{10}+z b_{01}+\sigma b_{20}+(\sigma z^2+2z c)b_{02})
\]
Pulling this back along $B\hookrightarrow B_I$ gives the subscheme of $U_\lambda\times B$
cut out by the ideal
\[
  (\alpha+\tau \beta+(\tau z^2+c^2)\gamma,
  \sigma \beta+(\sigma z^2+2z c)\gamma).
\]
One can verify $\sC^{[2]}$ restricts to an irreducible subvariety of
$U_\lambda\times B$ which has a 2-dimensional singular locus.
\begin{Verbatim}[frame=single,label=\textsc{Singular}]
> ring R = 0,(s,t,z,c,b,g),dp;
> poly p = sb + (sz2+2zc)*g;
> factorize(p);
[1]:
   _[1]=1
   _[2]=sz2g+2zcg+sb
[2]:
   1,1
> ideal I = p, jacob(p);
> dim(std(I));
2
\end{Verbatim}

\subsection{Family of cubics}
Consider now the family defined by
\[\sC=\spec\frac{\C[\alpha,\beta,\gamma,x,y]}{\alpha+\beta y+\gamma x^3}\]
over the same base $B$. The discriminant locus is $\set{\alpha=\beta=0}$,
where the fibre is the non-reduced $y$-axis $\spec(\C[x,y]/(x^3))$.

Take the index set $I=\set{(0,0),(1,0),(0,1),(2,0),(1,1),(0,2),(3,0)}$.
Now, let $\lambda$ be the partition $21$ of 3. One can check that
\[U_\lambda=\spec\C[a,b,c,d,e,f]\]
and the restriction of the universal family to the chart $U_\lambda$ is
cut out by the ideal
\begin{align*}
  (
  &x^2-ax-by+ad-bc+bf-d^2,\\
  &xy-cx-dy+cd-be,\\
  &y^2-ex-fy+ae-c^2+cf-de)
    \triangleleft\Os_{U_\lambda}\otimes_\C\C[x,y].
\end{align*}
Consider the morphism
\[
  \rho\colon\Os_{U_\lambda}
  \oplus\Os_{U_\lambda}x\oplus\Os_{U_\lambda}y
  \oplus\Os_{U_\lambda}x^2\oplus\Os_{U_\lambda}xy\oplus\Os_{U_\lambda}y^2
  \oplus\Os_{U_\lambda}x^3\to\pi_\ast\Os_{\U}\vert_{U_\lambda}
\]
Under the isomorphism
$\pi_\ast\Os_{\U}\vert_{U_\lambda}\cong\Os_{U_\lambda}\oplus\Os_{U_\lambda}x
\oplus\Os_{U_\lambda}y$, the matrix of this morphism is the transpose
of the following matrix.
\[\begin{pmatrix}
  1&0&0\\0&1&0\\0&0&1\\
  ad-bc+bf-d^2&-a&-b\\
  cd-be&-c&-d\\
  ae-c^2+cf-de&-e&-e\\
  ad^2-a^2d+abc-abf+b^2e-bcd&a^2-ad+2bc-bf+d^2&ab+bd
\end{pmatrix}\]
Thus, $\sC^{[3]}\cap(U_\lambda\times B)$ is cut out by the ideal
\begin{align*}
  (
  &\alpha+\gamma(ad^2-a^2d+abc-abf+b^2e-bcd),\\
  &\gamma(a^2-ad+2bc-bf+d^2),\\
  &\beta+\gamma(ab+bd))\triangleleft\C[a,b,c,d,e,f,\alpha,\beta,\gamma].
\end{align*}
Therefore,
\[
  \sC^{[3]}\cap(U_\lambda\times B)
  \cong\spec\frac{\C[a,b,c,d,e,f,\gamma]}
  {\left(\gamma(a^2-ad+2bc-bf+d^2)\right)}.
\]
The Hilbert scheme has two components, but the other component is the one
parametrising subschemes on the affine plane, which is the fibre of $\sC\to B$
above 0. After removing this degenerate fibre from $\sC\to B$,
one is left with an irreducible Hilbert scheme cut out by $a^2-ad+2bc-bf+d^2$.
This space has 3-dimensional singular locus:
\begin{Verbatim}[frame=single,label=\textsc{Singular}]
> ring R = 0,(a,b,c,d,e,f,g),dp;
> poly p = a2-ad-bf+d2;
> ideal I = p, jacob(p);
> dim(std(I));
3
\end{Verbatim}

\bibliographystyle{alpha}
\bibliography{rel-hilb-curves}

\end{document}